\documentclass[a4paper,11pt]{article}

\input{diagrams}

\newcommand{\op}{\ensuremath{\mbox{\hspace{1pt}{\scriptsize op}}}}

\newcommand{\cat}[1]{\ensuremath{\mbox{\bfseries {\upshape {#1}}}}}

\newcommand{\cl}[1]{\ensuremath{\mathcal {#1}}}
\newcommand{\bb}[1]{\ensuremath{\mathbb {#1}}}

\newcommand{\kml}{Kelly-Mac~Lane }

\newcommand{\tensordot}{\ensuremath{\otimes \cdots \otimes}}

\newcommand{\lra}{\ensuremath{\longrightarrow}}
\newcommand{\map}[1]{\ensuremath{\stackrel{{#1}}{\lra}}}

\newtheorem{theorem}{Theorem}[section]
\newtheorem{proposition}[theorem]{Proposition}
\newtheorem{corollary}[theorem]{Corollary}
\newtheorem{lemma}[theorem]{Lemma}

\newenvironment{prf}{\vspace{2ex}\begin{sloppypar}{\noindent\upshape
{\bfseries Proof. }}} {{\hspace*{\fill}
$\Box$}\end{sloppypar}\vspace{2ex}}

\newenvironment{prfof}[1]{\vspace{2ex}\begin{sloppypar}{\noindent
\upshape{\bfseries Proof of {#1}. }}} {{\hspace*{\fill}
$\Box$}\end{sloppypar}\vspace{2ex}}

\newcommand{\numroman}{\renewcommand{\labelenumi}{\roman{enumi})}}
\newcommand{\numarabic}{\renewcommand{\labelenumi}{\arabic{enumi})}}

\newcommand{\pica}{\begin{center} \input}
\newcommand{\picz}{\end{center}}
\newcommand{\length}[1]{\setlength{\unitlength}{#1}}

\newlength{\leng}
\newlength{\fontleng}
\newcommand{\sunit}{\setlength{\unitlength}{1mm}}

\sunit

\sunit

\sunit

\sunit

\sunit

\sunit

\sunit

\sunit

\sunit

\usepackage{amsfonts}
\usepackage{amssymb}
\usepackage{amsmath}
\usepackage{fancyheadings}
\usepackage{latexcad}

\begin{document}

\title{A relationship between trees and Kelly-Mac Lane graphs}
\author{Eugenia Cheng\\ \\Department of Pure Mathematics, University
of Cambridge\\E-mail: e.cheng@dpmms.cam.ac.uk}
\date{October 2002}
\maketitle

\begin{abstract}

We give a precise description of combed trees in terms of Kelly-Mac~Lane
graphs.  We show that any combed tree is uniquely expressed as an allowable
Kelly-Mac~Lane graph of a certain shape.  Conversely, we show that any such
Kelly-Mac~Lane graph uniquely defines a combed tree.  

\end{abstract}

\setcounter{tocdepth}{3}
\tableofcontents

\section*{Introduction}
\addcontentsline{toc}{section}{Introduction}

In this paper we show how trees may be expressed as allowable Kelly-Mac~Lane
graphs of a certain shape.  

The trees in question are those arising in \cite{bd1} and \cite{che7}
to express configurations for composing higher-dimensional cells in
the theory of opetopic $n$-categories.  Kelly-Mac~Lane graphs are introduced in
\cite{km1} to study coherence for symmetric monoidal closed
categories.  The main result of this paper states that a tree is {\em
precisely} an allowable Kelly-Mac~Lane graph of shape
	\[X_{m_1} \tensordot X_{m_k} \lra X_{(\sum\nolimits_i
	m_i-k+1)}\]
where each $m_i \geq 0$ and for any $m \geq 0$ we write
	\[X_m = [1^{\otimes m},1].\]
	
We begin, in Section~\ref{loopsap}, by giving a minimal account of
the theory of Kelly-Mac~Lane graphs, including no more than what is
required for the purposes of this paper.  We refer the reader
to \cite{km1} for the full details, noting that for the purposes of
this work we consider only the strict monoidal case.

We then recall the trees in question, as defined in \cite{che7}.  We
give an informal description in Section~\ref{treeback}, and in
Section~\ref{treeformal} the formal description that paves the way
for the ensuing results.  In Section~\ref{treegraph} we show how to
express a tree as a graph, not {\em a priori} allowable, and in
Section~\ref{treecomp} we characterise tree composition in this new
framework.  

This enables us to prove, in Section~\ref{every}, that the graph of a
tree is allowable.  Finally, in Section~\ref{all}, we show that every
allowable graph of the correct shape is a tree, giving the main result
of this work.

\subsection*{Further work}
\numarabic

\begin{enumerate}

\item The trees in question arise in the `slice' construction
described in \cite{bd1} and \cite{che7} in the construction of
opetopes and definition of opetopic $n$-category.  Thus the new way of
expressing trees as described in this paper can be used to give a
description of the category of opetopes; this is described in
\cite{che13}.  

\item Blute(\cite{blu1}) has established a relationship between Kelly-Mac~Lane
graphs and the proof nets of Linear Logic, so the material in the
present work should in turn give a relationship between opetopes and
proof nets.

\item Finally we note that the allowable Kelly-Mac~Lane graphs are the
morphisms of the free symmetric monoidal closed category on one
object.  Although we do not need to use this fact in this paper, it
should give a more abstract approach to this material, the
significance of which is currently unclear.  

\end{enumerate}

\bigskip {\bfseries Acknowledgements}

This work was supported by a PhD grant from EPSRC.  I would like to thank
Martin Hyland and Tom Leinster for their support and guidance.

\section{Background on Kelly-Mac~Lane Graphs} \label{loopsap}

In this section we give a brief account of the theory of
Kelly-Mac~Lane graphs.  Note that we will only be concerned with the
strict monoidal version for the purposes of this paper.  

In \cite{km1}, Kelly and Mac Lane study
coherence for symmetric monoidal closed categories.  In brief, a
symmetric monoidal closed category is a symmetric monoidal
category $\cl{C} = (\cl{C}, \otimes, I, a, b, c)$ equipped, in
addition, with a functor
    \[ [\ ,\ ]:\cl{C}^{\op}\times\cl{C} \lra \cl{C}\]
and natural transformations
    \[d = d_{AB}:A \lra [B, A \otimes B]\]
    \[e = e_{AB}:[A,B]\otimes A \lra B\]
satisfying certain axioms.  (Here $a$, $b$ and $c$ are the natural
isomorphisms for associativity, unit and symmetric action
respectively.) In particular we have a natural isomorphism
    \[\pi: \cl{C}(A \otimes B, C) \lra \cl{C} (A, [B,C]).\]
Kelly and Mac Lane refer to such categories simply as {\em closed
categories} and we do the same.

Kelly and Mac Lane introduce a notion of graph which enables a
partial solution to the question: when does a diagram in a closed
category commute?   In fact we are not concerned with the
coherence question here, so we only give the construction of the
graphs and state one theorem from \cite{km1} which will later be
useful.

Kelly and Mac Lane define a category $G$ whose objects are {\em
shapes} and whose morphisms are {\em graphs}; this is seen to be a
closed category.  They then define a subcategory whose morphisms
are the {\em allowable morphisms}.  These are defined as precisely
those morphisms of $G$ demanded by the symmetric monoidal closed
structure.

\subsection{Shapes}
\numarabic

We define {\em shapes} by the following inductive rules:
\begin{enumerate}
\item $I$ is a shape
\item 1 is a shape
\item if $S$ and $T$ are shapes then so is $S \otimes T$
\item if $S$ and $T$ are shapes then so is $[S,T]$
\end{enumerate}

Thus shapes are formal objects built from 1, $I$, $\otimes$ and
$[\ ,\ ]$.

We assign to each shape $T$ a {\em variable set} $v(T)$ which may
be considered as a list of $+$'s and $-$'s, defined inductively as
follows:

\begin{enumerate}
\item $v(I) = \emptyset$
\item $v(1) = \{+\}$
\item $v(T \otimes S) = v(T) \coprod v(S)$
\item $v([T,S]) = v(T)^{\op} \coprod v(S)$
\end{enumerate}

Here $\coprod$ is the concatenation of lists and $v(T)^{\op} $ is
$v(T)$ with all signs reversed.  Kelly and Mac Lane write
    \[v(T) \coprod v(S) = v(T) \hat{+} v(S)\]
    \[v(T)^{\op}  \coprod v(S) = v(T) \tilde{+} v(S)\]
and call these the {\em ordered sum} and {\em twisted sum}
respectively.  The sign of each variable is called its {\em
variance}.

In fact we only need the strict monoidal version of this theory.
That is, we put
    \[(T\otimes S) \otimes R = T \otimes (S \otimes R)\]
and
    \[T \otimes I = T.\]

For example,
     \[ [\ [1,1] \otimes 1\otimes 1 \ , \ I \ ] \otimes 1\]
is a shape with
    \[v(T) = \{+,-,-,-,+\}.\]

\subsection{Graphs} \label{graphs}

A {\em graph} $T \lra S$ is defined to be a fixed point free
pairing of the variables in $T$ and $S$ such that paired elements
have opposite variances in $v(T)^{\op}  \coprod v(S)$.  (Kelly and
Mac Lane refer to such paired elements as ``mates''.)
Equivalently, this is a bijection between the $+$'s and the $-$'s
in $v(T)^{\op}  \coprod v(S)$.

For example, the following is a graph
    \[ [ \ [1,1] \otimes 1 \otimes 1 \ , \ I \ ] \ \otimes \ 1
    \lra [ \ 1\otimes1 \ , \ 1\otimes [1,1] \ ]\]
showing variances: \pica pic52.lp .\picz

Graphs are composed in the obvious way, so that shapes and graphs
form a category $G$.  Moreover, $G$ has the structure of a closed
category as follows.  $\otimes$ and $[\ ,\ ]$ are defined on
graphs in the obvious way, and the constraints are given by the
following graphs:

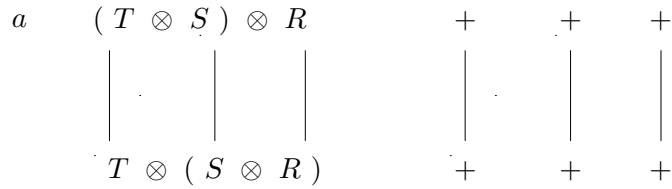
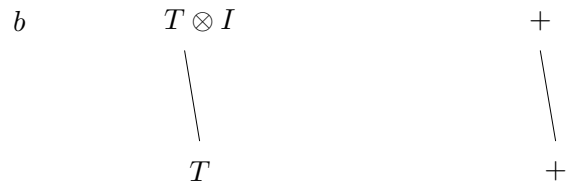
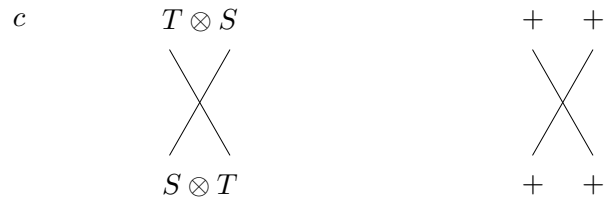
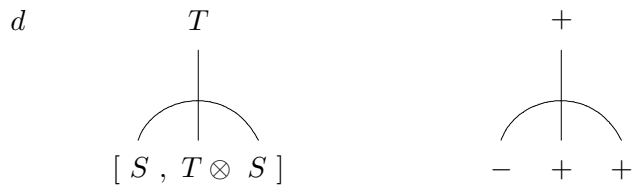
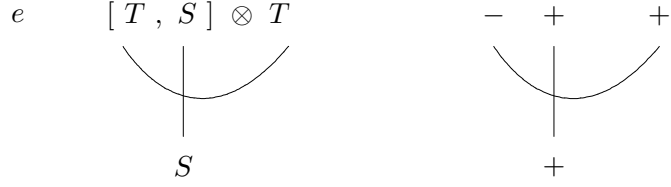
\input pic65b.lp \input pic70b.lp

\input pic66b.lp \input pic71b.lp

\input pic67b.lp \input pic72b.lp

\input pic68.lp
 \input pic73b.lp

\input pic69.lp \input pic74b.lp .

The diagrams on the right give variances, showing that these are
indeed graphs; note that in the twisted sum the variances of the
domain are reversed.  For the strict monoidal version we have
$a=1$ and $b=1$.

Observe that we realise \kml graphs as pictorial graphs by joining
paired objects up with an edge.  In the diagrams above, the objects
are in fact shapes, so the drawn edges in fact represent multiple
edges as necessary.

In fact there is a notion of graphs labelled in a category \bb{C}
(see \cite{che13}); these are the morphisms of a category which we
will call $K\bb{C}$.  Then the graphs above may be considered as
graphs labelled in the category \cat{1}.  So for consistency we write
$G = K\cat{1}$.

\subsection{Allowable morphisms}

The {\em allowable morphisms} are then defined to be the smallest
class of morphisms of $K\cat{1}$ satisfying the following
conditions:

\begin{enumerate}
\item For any $T,S,R$ each of the following morphisms is in the
class:
\[ \begin{array}{rcl}
    1 & : & T \lra T \\
    a & : & (T\otimes S) \otimes R \lra T \otimes (S \otimes R)\\
    a^{-1} & : & T \otimes (S \otimes R) \lra (T \otimes S) \otimes
    R\\
    b & : & T \otimes I \lra T\\
    b^{-1} & : & T \lra T \otimes I\\
    c & : & T \otimes S \lra S \otimes T.
\end{array}\]

\item For any $T, S$ each of the following morphisms is in the
class:
\[ \begin{array}{rcl}
    d & : & T \lra [S,T\otimes S]\\
    e & : & [T,S] \otimes T \lra S.
\end{array} \]

\item If $f: T\lra T'$ and $g:S\lra S'$ are in the class so is
    \[f\otimes g: T\otimes S \lra T' \otimes S'.\]

\item If $f:T \lra T'$ and $g:S \lra S'$ are in the class then so is
    \[[f,g] : [T',S] \lra [T,S'].\]

\item If $f: T \lra S$ and $g:S \lra R$ are in the class then so is
$gf: T \lra R$.

\end{enumerate}

We write $A\cat{1}$ for the category of shapes and allowable
morphisms.

The main theorem of \cite{km1} that we use is as follows:

\begin{theorem} \label{noloops}
If $f:T \lra S$ and $g:S \lra R \ \in G$ are allowable then they
are compatible, that is, composing them gives no closed loops.
\end{theorem}

For the proof (an induction over structure), see \cite{km1}.

\subsection{Duality} \label{duality}

Since $K\cat{1}$ is closed, given any graph
    \[\xi:S \otimes T \lra U \ \in K\cat{1}\]
there is a unique dual
    \[\bar{\xi}:S \lra [T,U]\]
so in particular, given a graph
    \[\alpha:S\lra T\]
there is a unique dual
    \[\bar{\alpha}: I \lra [S, T].\]
We will eventually be concerned with graphs of the form
    \[\alpha: A_1 \tensordot A_k \lra B;\]
it is sometimes convenient or indeed necessary to use the dual
    \[\bar{\alpha} : I \lra [\ A_1 \tensordot A_k \ , \ B \ ]\]
and we may refer to either of these graphs as $\alpha$ when the
exact form is not relevant.

\sunit
\section{Informal description of trees} \label{treeback}

We consider unlablled, `combed' trees, with ordered nodes.  For
example the following is a tree:

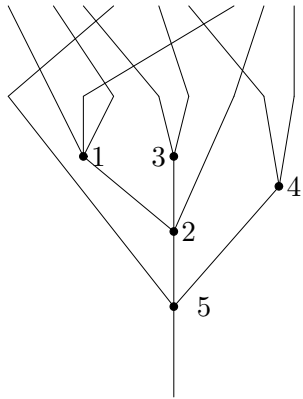
\begin{center}
\input pic01.lp .
\end{center}

\renewcommand{\labelenumi}{\roman{enumi})}

Explicitly, a tree $T=(T, \rho, \tau)$ consists of
\begin{enumerate}
\item A planar tree $T$
\item A permutation $\rho \in \cat{S}_l$ where $l=$ number of leaves of
$T$
\item A bijection $\tau: \{\mbox{nodes of $T$}\} \lra \{1, 2,
\ldots, k\}$ where $k=$ number of nodes of $T$; equivalently an
ordering on the nodes of $T$.
\end{enumerate}


Note that there is a `null tree' with no nodes
\begin{center}
\input pic53.lp .
\end{center}

%

\section{Formal description of trees}
\label{treeformal}

In this section we give a formal description of the above trees,
characterising them as connected
graphs with no closed loops (in the conventional sense of `graph'). 
This enables us, in Section~\ref{treegraph}, to express a tree as a
morphism in $K\cat{1}$; it also enables us, in Section~\ref{all}, to
show that all {\em allowable} graphs of the correct shape arise in
this way.

We consider a tree with $k$ nodes $N_1, \ldots, N_k$ where $N_i$
has $m_i$ inputs and one output. Let $N$ be a node with
$(\sum\limits_i m_i )- k+1$ inputs; $N$ will be used to represent
the leaves and root of the tree.

Then a tree is given by a bijection
    \[\coprod_i\{\mbox{inputs of $N_i$}\} \coprod \{\mbox{output of $N$}\} \lra
    \coprod_i\{\mbox{output of $N_i$}\} \coprod \{\mbox{inputs of
    $N$}\}\]
since each input of a node is either connected to a unique output
of another node, or it is a leaf, that is, input of $N$. Similarly
each output of a node is either attached to an input of another
node, or it is the root, that is, output of $N$.

We express this formally as follows.

\begin{lemma} \label{tree1}
Let $T$ be a tree with nodes $N_1, \ldots, N_k$, where $N_i$ has
inputs $\{x_{i1}, \ldots, x_{im_i}\}$ and output $x_i$.  Let $N$
be a node with inputs $\{z_1, \ldots, z_l\}$ and output $z$, with
    \[l=(\sum\limits_{i=1}^k m_i) - k +1.\]
    Then $T$ is given by a
bijection
    \[\alpha \ : \ \coprod_i\{x_{i1},\ldots, x_{im_i}\} \coprod \{z\} \lra
    \coprod_i\{x_i\} \coprod \{z_1, \ldots, z_l\}.\]
\end{lemma}

\begin{prf}  We construct the bijection $\alpha$.

Consider $x_{ij}$ on the left hand side.  This is the $j$th input
of $N_i$, which is either
\begin{enumerate}
\item joined to the output of a unique $N_r$, in which case
$\alpha(x_{ij})=x_r$, or
\item the $p$th leaf of the tree, in which case
$\alpha(x_{ij})=z_p$.
\end{enumerate}  Finally, $z$ is the root of the tree, so must be
the output of a unique $N_r$, so $\alpha(z)=x_r$.

For the inverse, consider $x_r$ on the right hand side.  This is
the output of the $r$th node, so is either
\begin{enumerate}
\item joined to the $j$th input of a unique $N_i$, in which case
$\alpha^{-1}(x_r) = \alpha(x_{ij})$, or
\item is the root of the tree, in which case $\alpha^{-1}(x_r)=z$.
\end{enumerate}

Each $z_r$ is a leaf of the tree, so must be the $j$th input of a
unique $N_i$, so $\alpha^{-1}(z_r)=x_{ij}$.

$\alpha^{-1}$ thus defined is inverse to $\alpha$, so $\alpha$ is
a bijection.

Note that if $k=0$ we have the null tree with no nodes; then $l=1$
and $N$ has one input $z_1$.  Then the bijection $\alpha$ is given
by $\alpha(z)=z_1$.

\end{prf}

For example, consider \pica pic21.lp \input
pic22.lp \input pic23.lp .\picz  Then a tree \pica
pic24.lp \picz is given by the following bijection:

\[\begin{array}{rcl}
    \{x_1, x_2, x_3, y_1, y_2, z\} & \lra & \{x,y,z_1, z_2, z_3,
    z_4\}\\
    x_1 & \longmapsto & z_1 \\
    x_2 & \longmapsto & z_3 \\
    x_3 & \longmapsto & z_4 \\
    y_1 & \longmapsto & x \\
    y_2 & \longmapsto & z_2 \\
    z & \longmapsto & y .\\
    \end{array}\]

For the converse, every such bijection gives a graph, but it is
not necessarily a tree.  For example

\[\begin{array}{rcl}
    x_1 & \longmapsto & y \\
    x_2 & \longmapsto & z_3 \\
    x_3 & \longmapsto & z_4 \\
    y_1 & \longmapsto & x \\
    y_2 & \longmapsto & z_2 \\
    z & \longmapsto & z_1 \\
    \end{array}\]
gives the following graph: \pica pic25.lp .\picz  So we need to
ensure that the resulting graph has no closed loops; the use of
the `formal' node $N$ then ensures connectedness.  We express this
formally as follows.

\begin{lemma} \label{loopspropb}
Let $N_1, \ldots, N_k, N$ be nodes where $N_i$ has inputs
$\{x_{i1}, \ldots, x_{im_i}\}$ and output $x_i$, and $N$ has
inputs $\{z_1, \ldots, z_l\}$ and output $z$, with
$l=(\sum\limits_{i=1}^k m_i) - k +1$. Let $\alpha$ be a bijection
    \[\coprod_i\{x_{i1},\ldots, x_{im_i}\} \coprod \{z\} \lra
    \coprod_i\{x_i\} \coprod \{z_1, \ldots, z_l\}.\]
Then $\alpha$ defines a graph with nodes $N_1, \ldots, N_k$.
\end{lemma}

\begin{lemma} \label{loopspropb2}
Let $\alpha$ be a graph as above.  Then $\alpha$ has a closed loop
if and only if there is a non-empty sequence of indices
    \[\{t_1 , \ldots, t_n \} \subseteq \{1, \ldots, k\}\]
such that for each $2 \leq j \leq n$
    \[\alpha( x_{t_jb_j})=x_{t_{j-1}} \]
for some $1 \leq b_j \leq m_j$, and
    \[\alpha(x_{t_1b_1})= x_{t_n}\]
for some $1 \leq b_1 \leq m_1$.  \end{lemma}

\begin{prf} A closed loop in $\alpha$ is a sequence of nodes
    \[\{N_{t_1}, \ldots, N_{t_n}\}\]
such that for each $2 \leq j \leq n$, $N_{t_j}$ is joined to
$N_{t_{j-1}}$, and also $N_{t_1}$ is joined to $N_{t_n}$.

That is, for each $2 \leq j \leq n$, some leaf of $N_{t_j}$ is
joined to the root of $N_{t_{j-1}}$, and also some leaf of
$N_{t_1}$ is joined to $N_{t_n}$. This is precisely the case
described formally in the Lemma, with the $b_j$ giving the leaves
in question.
\end{prf}

For example in the above case we have
\[\begin{array}{rrcl}
    \alpha : &x_{11} & \longmapsto & x_2 \\
    &x_{12} & \longmapsto & z_3 \\
    &x_{13} & \longmapsto & z_4 \\
    &x_{21} & \longmapsto & x_1 \\
    &x_{22} & \longmapsto & z_2 \\
    &z & \longmapsto & z_1 \\
    \end{array}\]
which has a loop given by indices $\{1,2\}$, since
    \[ \alpha(x_{21}) = x_1 \mbox{\ \ and \ \ } \alpha(x_{11}) = x_2.\]

Note that a graph with no nodes cannot satisfy the above condition
since the sequence $\{N_{t_1}, \ldots, N_{t_n}\}$ is required to
be non-empty.

\begin{corollary} \label{loopspropb3}
A tree with nodes $N_1, \ldots, N_k$ is precisely a bijection
$\alpha$ as in Lemma~\ref{loopspropb}, such that there is no
sequence of indices as in Lemma~\ref{loopspropb2}.
\end{corollary}

\begin{prf} $\alpha$ defines a graph; this is a tree if and only if there is
no closed loop. Note that if $k=0$ we have a bijection
    \[\alpha: \{z\} \lra \{z_1\}\]
that is, the null tree.  \end{prf}

\section{Trees as morphisms in $K\cat{1}$}
\label{treegraph} We now show how trees may be expressed as
graphs.  Here we consider unlabelled trees; the labelled version
follows easily.

Let \cat{1} be the category with just one object and one
(identity) morphism.  We write the single object of \cat{1} as
$1$.  Then we express a node of a tree as the following object in
$K\cat{1}$
    \[X_m = [1 \otimes \ldots \otimes 1, 1] = [1^{\otimes m}, 1]\]
where $m$ is the number of inputs of the node.

%

Now consider a tree $T$ with (ordered) nodes $N_1, \ldots N_k$
where $N_i$ has $m_i$ inputs.  We show that this tree may be
represented as a morphism
    \[X_{m_1} \otimes \ldots \otimes X_{m_k} \map{\xi_T} X_l \ \in \
    K\cat{1}\]
using the formal description of trees as in
Section~\ref{treeformal}.

\begin{lemma} \label{tree2}
Let $T$ be a tree with $N_1, \ldots, N_k$ be nodes where $N_i$ has
inputs $\{x_{i1}, \ldots, x_{im_i}\}$ and output $x_i$.  Then $T$
is given by a morphism
    \[\xi_T : X_{m_1} \otimes \ldots \otimes X_{m_k} \lra X_l \ \in \
    K\cat{1}\]
where $l=(\sum\limits_{i=1}^k m_i) - k +1$. Note that if $k=0$
then the left hand side of the above expression becomes $I$.
\end{lemma}

\begin{prf}
Recall that a graph $\xi_T$ as above is precisely a bijection from
the $-$'s to the $+$'s in the twisted sum
    \[v(X_{m_1} \otimes \ldots \otimes X_{m_k}) \tilde{+} v(X_l).\]
By Lemma~\ref{tree1}, $T$ is given by a bijection
    \[\coprod_i\{x_{i1},\ldots, x_{im_i}\} \coprod \{z\} \lra
    \coprod_i\{x_i\} \coprod \{z_1, \ldots, z_l\}.\]
Observe that the elements of the left hand side of this expression
are precisely the $-$'s in the twisted sum above, and those of the
right hand side are precisely the $+$'s. \end{prf}

As in Section~\ref{treeformal}, the idea is that a tree is
constructed by identifying each node output with the node input to
which it is joined, unless it is the root; similarly each input is
identified with a node output unless it is a leaf.  This
identification gives the mates in the graph $\xi_T$, where the
codomain $X_l$ is representing the leaves and the root of the tree
$T$.

For example the following tree as described in
Section~\ref{treeformal}
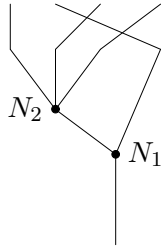
\begin{center}
\input pic02.lp
\end{center}
is be expressed as the following morphism in $K\cat{1}$
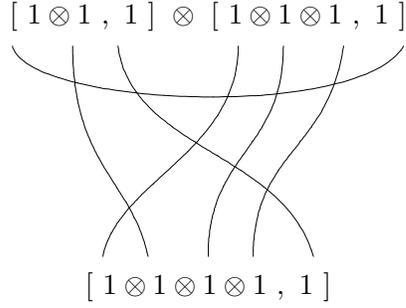
\begin{center}
\input pic03.lp
\end{center}
and the following representation giving variances shows that this
is indeed a graph:
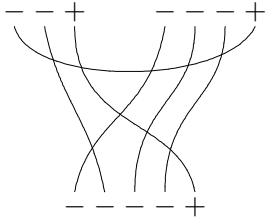
\begin{center}
\input pic04.lp
\end{center}
Formally, the graph for a tree $T$ as above is given as follows.
We write
    \[X_{m_i} = [\ A_{i1} \otimes \ldots \otimes
    A_{im_i}\ ,\  A_i\ ] \]
    \[X_l=[\ B_1, \otimes \ldots \otimes B_l, \
    , \ B\ ]\]
where each $A_{ij}, A_i, B_i, B = 1$ and in the twisted sum we
have variances
    \[\begin{array}{rl}
    v(A_{ij}) = +, & v(A_i) = -\\
    v(B_p) = -, & v(B) = +.
    \end{array}\]
Then the graph $\xi_T$ is given as follows.

\begin{itemize}
\item considering node inputs
\end{itemize}
For each $i,j$, either
\begin{enumerate}
\item the $j$th input of $N_i$ is joined to the output of
$N_r$, say, in which case $A_{ij}$ is the mate of $A_r$, or
\item the $j$th input of $N_i$ is the $p$th leaf of the tree $T$,
in which case $A_{ij}$ is the mate of $B_p$ in $\xi_T$.
\end{enumerate}

\begin{itemize}
\item considering node outputs
\end{itemize}
For each $r$, either
\begin{enumerate}
\item the output of $N_r$ is the root of the tree, in which case $B_r$ is the
mate of $B$, or
\item the output of $N_r$ is joined to the $j$th input of $N_i$,
say, in which case $A_r$ is the mate of $A_{ij}$.
\end{enumerate}

Note that the null tree \pica pic53.lp \picz is a graph $I
\map{\xi} X_1$ as follows: \pica pic55.lp .\picz

So we have shown that every tree is given by a graph in
$K\cat{1}$; in Section~\ref{every} we show that any such graph is
allowable.  The proof is by induction, and the following section
enables us to makes the induction step.

\section{Composition of trees} \label{treecomp}

We now discuss two ways of composing trees:

\numarabic
\begin{enumerate}
\item {\em leaf-root composition} in which a leaf of one tree is
attached to the root of another, for example \pica pic05.lp \picz
\item {\em node-replacement composition} in which a node of one
tree is replaced by another tree, for example \pica pic06.lp \picz
\end{enumerate}

In the first case the inputs of the tree are considered to be the
leaves, and the output the root; note that an issue of
node-ordering arises, so that this `composition' is not
associative.  However, it facilitates the induction argument in
Section~\ref{every}, which is why we discuss it here.

In the second case the inputs are the nodes, and the output a node
with one input edge for each leaf of the tree.  This form of
composition is used in later work (\cite{che13}) when we describe the
construction of opetopes in the framework of Kelly-Mac~Lane graphs.  

We show how each of these forms of composition arises for trees
represented as graphs as in Section~\ref{treegraph}

Recall that a tree is expressed as a morphism
    \[X_{m_1}\otimes \cdots \otimes X_{m_k} \lra X_l \ \in K\cat{1}.\]
Now in general, given any morphisms in $K\cat{1}$
    \[B_1 \tensordot B_n \map{f} A_p\]
    \[A_1 \tensordot A_m \map{g} A\]
for some $1 \leq p \leq m$, we may form the composite
    \[f\circ(1 \tensordot 1 \otimes g \otimes 1 \tensordot 1)\]
which we write as
    \[g\circ_p f:A_1 \tensordot A_{j-1} \otimes B_1 \tensordot B_n
    \otimes A_{p+1} \tensordot A_m \lra A.\]
Note that if $p$ is evident from the context we simply write
$g\circ f$.  

This composition gives node-replacement composition of trees.
Consider trees $S, T$ with graphs
    \[\xi_S: X_{s_1} \tensordot X_{s_n} \lra X_l\]
    \[\xi_T: X_{t_1} \tensordot X_{t_m} \lra X_k\ .\]
Then S may be composed at the $p$th node of $T$ if the number of
leaves of $S$ equals the number of inputs of the $p$th node, that
is, if $X_l=X_{t_p}$.  Then the graph for the composite tree is
given by
    \[\xi_S \circ_p \xi_T.\]

For example as above, suppose we have $p=2$ and \pica pic07.lp
\picz then we express this with graphs as follows \pica pic08.lp
\picz

In fact, considering the dual forms $\bar{\xi}_S$ and
$\bar{\xi}_T$, we see that this composite may also be expressed by
means of a `composition graph' $\xi$ as follows.  We have
    \[\bar{\xi}_S: I \lra [X_{s_1} \tensordot X_{s_n} \ , \ X_l]\]
    \[\bar{\xi}_T: I \lra [X_{t_1} \tensordot X_{t_m} \ , \ X_k].\]
Then $\xi$ is a graph
\[\begin{array}{c}[X_{t_1} \tensordot X_{t_m}\ , \ X_k]
    \otimes [X_{s_1} \tensordot X_{s_n} \ , \ X_l] \\
    \downarrow \\
    \mbox{\hspace{0mm}}[X_{t_1} \tensordot X_{t_{p-1}} \otimes X_{s_1}
    \tensordot X_{s_n} \otimes X_{t_{p+1}} \tensordot X_{t_m} \ , \ X_k]
\end{array}\]
where $X_l$ is joined to $X_{t_p}$ in the domain, and for all
other $j$, $X_j$ in the domain is joined to $X_j$ in the codomain.

\bigskip

We now consider leaf-root composition.  Consider trees $S,T$ as
above.  We seek to attach the root of $S$ to the $q$th leaf of
$T$, and we adopt the convention that the nodes of $S$ are then
listed before those of $T$ in the final tree.

This is achieved in $K\cat{1}$ by placing the graphs $\xi_S$ and
$\xi_T$ side by side, that is, taking their tensor product, and
composing the result with a `composition graph' that joins up the
correct leaf and root as required.  We write
    \[X_l=[A_1 \tensordot A_l, A]\]
    \[X_k=[B_1 \tensordot B_k, B]\]
    \[X_{l+k-1}=[C_1 \tensordot C_{l+k-1},C]\]
and the `composition graph' as
    \[\xi:X_l \otimes X_k \lra X_{l+k-1}.\]
The idea is that the leaves of $S$ are inserted into the list of
leaves of $T$ at the $q$th place to give
    \[[B_1 \tensordot B_{q-1} \otimes A_1 \tensordot A_l \otimes
    B_{q+1} \tensordot B_k\ ,\ B]\]
so the composition graph $\xi$ is given as follows: \numroman

\begin{enumerate}
\item the mate of $A$ is $B_q$
\item the mate of $B$ is $C$
\item for $1\leq i \leq l$ the mate of $A_i$ is $C_{q+i-1}$
\item for $1\leq i \leq q-1$ the mate of $B_i$ is $C_i$
\item for $q+1 \leq i \leq k$ the mate of $B_i$ is $C_{l+i-1}$.
\end{enumerate}

For example, suppose we have $q=2$ with \pica pic10.lp \picz then
this is represented by the following graph in $K\cat{1}$: \pica
pic11.lp \input pic12.lp \picz

Note that we could adopt a different convention for ordering the
nodes of the composite tree, using $\xi_T \otimes \xi_S$.  Of
course, neither convention yields an associative composition, but
since we are not at this time trying to form a category (or
multicategory) of such trees, we do not pursue this matter here.

\section{The graph of a tree is allowable}
\label{every}

We have shown how any tree is represented by a graph.  We now show
that any such graph is {\em allowable}.

\begin{proposition} \label{loopspropa}
Given a tree $T$ as above, the graph $\xi_T$ is allowable.
\end{proposition}

\begin{prf}
By induction on the height of trees.  Here the height of a tree is
the maximum number of nodes on any path from a leaf to the root. A
tree of height 0 is the null tree \pica pic53.lp \picz represented
by the graph \pica pic13.lp \picz which is the morphism
    \[I \map{d_{I1}} [1, I\otimes 1] = [1,1]\]
which is allowable.

A tree of height 1 is just a node \pica pic14.lp \picz which is
represented by an identity graph \pica pic15.lp \picz which is
allowable.

A tree of height $h \ge 1$ may be considered as a composite
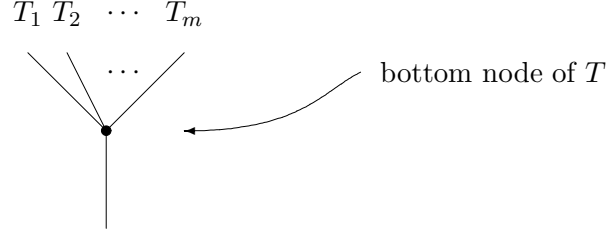
\length{1.3mm} \pica pic16.lp \picz where the $T_i$ are subtrees
of $T$; by construction they have height $\le h$.  So by induction
each of these is represented by an allowable graph. \length{0.9mm}

It is therefore sufficient to show that leaf-root composition of
allowable graphs gives an allowable graph.  Note that leaf-root
composition as defined in Section~\ref{treecomp} will not
necessarily give the correct node ordering on the final tree;
however, this can be achieved by composing with symmetries as
necessary.  This will not affect the allowability of the graph
since symmetries are allowable graphs, and composites of allowable
graphs are allowable.

Furthermore, since tensors and composites of allowable graphs are
allowable, it is sufficient to show that all `composition graphs'
$\xi$ as defined in Section~\ref{treecomp} are allowable.

Since any permutation may be written as a composite of
transpositions, and is therefore allowable, we may assume without
loss of generality that $q=1$ in the composition.  So it is
sufficient to show that any graph $\xi$ of the following form is
allowable.  Writing
    \[X_{m_1}=[A_1 \tensordot A_{m_1}, A]\]
    \[X_{m_2}=[B_1 \tensordot B_{m_2}, B]\]
    \[X_{m_1 + m_2 -1}=[C_1 \tensordot C_{m_1 + m_2 -1},C]\]
then
    \[\xi:X_{m_1} \otimes X_{m_2} \lra X_{m_1 + m_2 -1}\]
is given as follows. \numroman

\begin{enumerate}
\item the mate of $A$ is $B_1$
\item the mate of $B$ is $C$
\item for all $1\leq i \leq m_1$ the mate of $A_i$ is $C_i$
\item for $2 \leq i \leq m_2$ the mate of $B_i$ is $C_{m_1 + i -1}$.
\end{enumerate}

So $\xi$ has the form \pica pic17.lp \picz

Writing
    \[A_1 \tensordot A_{m_1} = \bar{A}\]
    \[B_2 \tensordot B_{m_2} = \bar{B}\]
we may abbreviate this as \pica pic18.lp \picz which may be
written as the following composite of allowable graphs: \pica
pic19.lp
\input pic20.lp \picz so $\xi$ is allowable as required.
\end{prf}

\section{Every allowable graph is a tree} \label{all}

We have seen that every tree is represented by a unique graph, and
that this graph is allowable. In this section we prove the
converse, that every allowable graph of the correct shape
represents a unique tree.

We now use the characterisation of trees as in
Section~\ref{treeformal}.  As in that section, for the converse we
see that every morphism
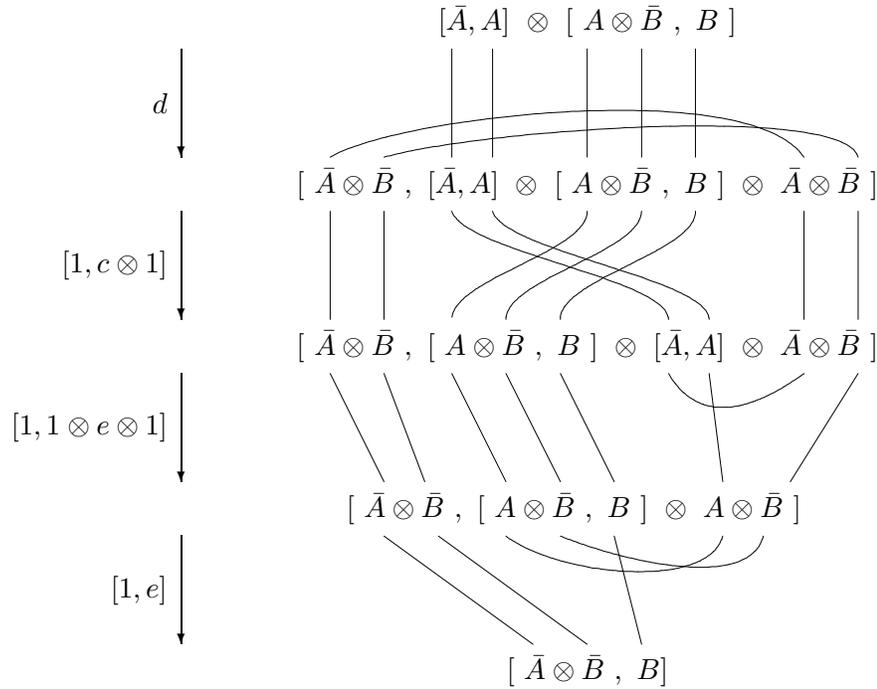
    \[X_{m_1}\otimes \cdots \otimes X_{m_k} \lra X_l \ \in K\cat{1}\]
gives a graph but that it is not necessarily a tree; we need to
ensure that the resulting graph has no closed loops.  We copy
Lemmas~\ref{loopspropb} and \ref{loopspropb2}, ``translating''
them into the language of closed categories.  Note that the word
`graph' is used in the ordinary sense; for clarity we refer to
Kelly-Mac~Lane graphs as `morphisms in $K\cat{1}$'.

\begin{lemma} \label{tree3}
Let $N_1, \ldots, N_k$ be nodes where $N_i$ has inputs
    \[\{A_{i1}, \ldots, A_{im_i}\}\]
and output $x_i$.  Let $\xi$ be a morphism
    \[\xi:X_{m_1} \otimes \ldots \otimes X_{m_k} \lra X_l \ \in \
    K\cat{1}\]
where $l=(\sum\limits_{i=1}^k m_i) - k +1$. Then $\xi$ defines a
graph with nodes $N_1, \ldots, N_k$. \end{lemma}

\begin{lemma} \label{tree4}
Let $\xi$ be a graph as above.  Then $\xi$ has a closed loop if
and only if there is a non-empty set of indices
    \[\{t_1 , \ldots, t_n \} \subseteq \{1,\ldots, k\}\]
such that for each $2 \leq j\leq n$ the mate of $A_{t_{j-1}}$
under $\xi$ is $A_{t_jb_j}$ and the mate of $A_{t_n}$ is
$A_{t_1b_1}$ for some $1 \leq b_j \leq m_j$. \end{lemma}

\begin{proposition} \label{loopspropd}
If there is a set of indices $\{t_1, \ldots t_n\}$ as above then
$\xi$ is not allowable.  \end{proposition}

\begin{corollary} Let $\xi$ be a morphism as above.  Then $\xi$ is
a tree if and only if it is allowable. \end{corollary}

To prove this we will use Theorem~\ref{noloops} (Theorem~2.2 of
\cite{km1}) which states that if two composable morphisms are
allowable then they are compatible, that is, composing them does
not result in any closed loops.  So to show that $\xi$ as above is
not allowable, we aim to construct an allowable morphism $\eta$
such that $\eta$ and $\xi$ are not compatible.  The following
lemma provides us with such a morphism.

\begin{lemma} \label{loopsprope}
Write $X_k = [A_1 \tensordot A_k, A]$ with $A_i, A=1$ and let $1
\leq p \leq k$.

Then there is an allowable morphism
    \[ \theta_p: [A_1 \tensordot A_{p-1} \otimes A_{p+1} \tensordot A_k, I]
    \lra X_k\]
with graph \pica pic26.lp .\picz \end{lemma}

\begin{prf} Write $Y=A_1 \tensordot A_{p-1} \otimes A_{p+1} \tensordot A_k$.
Since symmetries are allowable, it is sufficient to exhibit an
allowable morphism
    \[[Y,I] \lra [Y\otimes 1, 1]\]
with underlying graph \pica pic27.lp .\picz  We have the following
composite of allowable morphisms: \pica pic28.lp \picz which has
the underlying graph as required; since composites of allowable
morphisms are allowable, the composite is allowable.  \end{prf}

\begin{prfof}{Proposition~\ref{loopspropd}} To show that
    \[\xi: X_{m_1} \tensordot X_{m_k} \lra X_l\]
is not allowable we construct an allowable morphism
    \[\eta: T \lra X_{m_1} \tensordot X_{m_k}\]
such that $\eta$ and $\xi$ are not compatible, that is, composing
them produces a closed loop.

We aim to construct $\eta$ in such a way that for each $1\leq j
\leq n$ the mate of $A_{t_j}$ is $A_{t_jb_j}$ so that in the
composite graph we have the following closed loop: \pica pic29.lp
.\picz  We use morphisms of the form $\theta_p$ as given in
Lemma~\ref{loopsprope}.

Put $T=Y_1 \tensordot Y_k$ where
\[Y_i = [A_{t_j1}, \tensordot A_{t_j(b_j-1)} \otimes
    A_{t_j(b_j+1)} \tensordot A_{t_jm_i}\ ,\ I]\]
if $i=t_j$ for some $1\leq j \leq n$, and
    \[Y_i = X_{m_i}\]

We define $\eta$ as a tensor product
    \[f_1 \tensordot f_k\ :\ Y_1 \tensordot Y_k\  \lra\ X_{m_1}
    \tensordot X_{m_k}\]
where
    \[f_i = \left\{ \begin{array} {l@{\extracolsep{3em}}l}
    \theta_{b_j} & \mbox{if $i=t_j$ for some $1\leq j \leq n$}\\
    1 & \mbox{otherwise} \end{array} \right. \]
By Lemma~\ref{loopsprope} each $f_i$ is allowable, so $\eta$ is
allowable.

Since the mate of $A_{t_j}$ under $\theta_{b_j}$ is $A_{t_jb_j}$
we have a closed loop as  above, so $\eta$ and $\xi$ are not
compatible.  Since $\eta$ is allowable, it follows from
Theorem~\ref{noloops} that $\xi$ is not allowable.
\end{prfof}

Finally we sum up the results of this section in the following
proposition.

\begin{proposition} \label{loopspropg}
A tree is  a unique morphism of the form
    \[X_{m_1} \tensordot X_{m_k} \lra X_l \in K\cat{1}\]
and this morphism is allowable.  Conversely, any such allowable
morphism represents a unique tree.
\end{proposition}

%

\begin{corollary} \label{loopsproph}
A tree is a unique allowable morphism of the form
    \[I \lra [X_{m_1} \tensordot X_{m_1}, X_l] \in K\cat{1}.\]
Conversely, any such allowable morphism represents a unique tree.
\end{corollary}

\begin{prf} Follows from the closed structure of $K\cat{1}$.
\end{prf}

In order to make Proposition~\ref{loopspropg} and
Corollary~\ref{loopsproph} more precise, we seek an equivalence
between a `category of trees' and a `category of allowable
morphisms'.  In fact, trees of this form arise naturally by
considering configurations for composing arrows of a symmetric
multicategory.  That is, they arise from the `slicing' process as
defined in \cite{bd1} and \cite{che7}; the trees then appear
as arrows of the multicategory $I^{2+}$, and so as objects of
$I^{3+}$, forming a category $\bb{C}_3$.

So we may consider the slice construction using the
representation in closed categories.  In considering this for
constructing trees, we in fact deal with all the machinery used in
constructing $k$-opetopes for all $k\geq 0$, since these are
formed by iterating the construction.  This is the subject of
\cite{che13}.

\addcontentsline{toc}{section}{References}
\bibliography{bib0209}

\begin{thebibliography}{10}

\bibitem{bd1}
John Baez and James Dolan.
\newblock Higher-dimensional algebra {I}{I}{I}: $n$-categories and the algebra
  of opetopes.
\newblock {\em Adv. Math.}, 135(2):145--206, 1998.
\newblock Also available via {\tt http://math.ucr.edu/home/baez}.

\bibitem{blu1}
R.~Blute.
\newblock Linear logic, coherence and dinaturality.
\newblock {\em Theoretical Computer Science}, 115:3--41, 1993.

\bibitem{che13}
Eugenia Cheng.
\newblock The theory of opetopes via {K}elly-{M}ac~{L}ane graphs, October 2002.
\newblock E-print {\tt math.CT/0304288}.


\bibitem{che7}
Eugenia Cheng.
\newblock Weak $n$-categories: opetopic and multitopic foundations, October
  2002.
\newblock E-print {\tt math.CT/0304277}.


\bibitem{gir1}
Jean-Yves Girard.
\newblock Linear logic: its syntax and semantics.
\newblock In Laurent~Regnier, Jean-Yves~Girard, Yves~Lafont, editors, {\em
  Advances in Linear Logic}, volume 222 of {\em LMS}, pages 1--42. Cambridge
  University Press, 1995.

\bibitem{js1}
Andr{\'e} Joyal and Ross Street.
\newblock Braided tensor categories.
\newblock {\em Advances in Mathematics}, 102:20--78, 1983.

\bibitem{kl1}
G.~Kelly and M.~Laplaza.
\newblock Coherence for compact closed categories.
\newblock {\em Journal of Pure and Applied Algebra}, 19:193--213, 1980.

\bibitem{km1}
G.~M. Kelly and S.~{Mac Lane}.
\newblock Coherence in closed categories.
\newblock {\em Journal of Pure and Applied Algebra}, 1(1):97--140, 1971.

\bibitem{laf1}
Yves Lafont.
\newblock From proof-nets to interaction nets.
\newblock In Laurent~Regnier Jean-Yves~Girard, Yves~Lafont, editor, {\em
  Advances in Linear Logic}, volume 222 of {\em LMS}, pages 225--247. Cambridge
  University Press, 1995.

\bibitem{mel2}
Paul-Andr\'{e} Melli\`{e}s.
\newblock On double categories and multiplicative linear logic, June 1999.
\newblock to appear in Mathematical Structures in Computer Science.

\bibitem{shu1}
Mei~Chee Shum.
\newblock {\em Tortile Tensor Categories}.
\newblock PhD thesis, Macquarie University, New South Wales, November 1989.

\end{thebibliography}

\nocite{kl1}

\nocite{js1} 

\nocite{shu1}

\nocite{mel2}

\nocite{laf1}

\nocite{gir1}

\end{document}